\documentclass[12pt]{cimsart}               
\received{Month 200X}       
\yearofpublication{200X}                \volume{000}
\startingpage{1}                        \cccline{00}

\authorheadline{V. Bangert, M. Katz}
\titleheadline{Stable systolic inequalities}

 \usepackage{amssymb}                         
 \usepackage{amsopn}                         

\newtheorem{lemma}[theorem]{Lemma}
\newtheorem{corollary}[theorem]{Corollary}

\newtheorem{definition}[theorem]{Definition}
\newtheorem{question}[theorem]{Question}
\newtheorem{remark}[theorem]{Remark}
\newtheorem{example}[theorem]{Example}

\newcommand{\nc}{\newcommand}
\nc{\on}{\operatorname}
\def\tag#1#2{\hbox to\textwidth{#1\hfil$\displaystyle #2$\hfil}}
\nc{\df}{\on{\it df}}
\nc{\vol}{\on{vol}}
\nc{\sys}{\on{sys}}
\nc{\stsys}{\on{stsys}}
\nc{\conf}{\on{conf}}
\nc{\Imag}{\on{Im}}
\nc{\Hom}{\on{Hom}}
\nc{\PD}{\on{PD}}
\nc{\TM}{\on{TM}}
\nc{\rk}{\on{rk}}
\nc{\spt}{\on{spt}}
\nc{\norm}[1]{\| #1 \|}
\nc{\parallelleer}{\norm{\ }} 
\nc{\parallelh}{\norm h} 
\nc{\parallelk}{\norm k} 
\nc{\parallelx}{\norm x} 
\nc{\parallelhrr}{\norm {h_\RR}} 
\nc{\parallelom}{\norm \omega} 
\nc{\parallelomij}{\norm {\omega_{i_j}}} 
\nc{\parallelomx}{\norm {\omega_{x}}} 
\nc{\parallelpi}{\norm \pi} 
\nc{\parallelalf}{\norm \alpha} 
\nc{\parallelalfs}{\norm {\alpha_s}} 
\nc{\parallelalfi}{\norm {\alpha_i}} 
\nc{\parallelalfij}{\norm {\alpha_{i_j}}} 
\nc{\parallelbeta}{\norm \beta} 
\nc{\parallelbetat}{\norm {\beta_t}} 
\nc{\parallelhcapalf}{\norm {h \cap \alpha}} 
\nc{\parallelPDralf}{\norm {\PD_\RR(\alpha)}} 
\nc{\strichleer}{| \  |}
\nc{\NN}{\mathbb N}
\nc{\RR}{\mathbb R}
\nc{\rr}{\mbox{$\scriptstyle\mathbb R$}}
\nc{\ZZ}{\mathbb Z}
\newenvironment{otheorem}[1]{\par\bigskip\noindent\textsc{#1}\bgroup\it\ }%
                            {\egroup\par\bigskip}

\begin{document}                        



\title{Stable systolic inequalities and cohomology products}              


\author{Victor Bangert \\            
\affil Mathematisches Institut\\ Universit\"at Freiburg, Eckerstr.~1,79104 Freiburg,Germany\\  
 \\ AND \\
 \\Mikhail Katz\thanks{Supported by The Israel
Science Foundation (grant no.\ 620/00-10.0).  Partially supported by
the Emmy Noether Research Institute and the Minerva Foundation of
Germany.
}\\               
 \affil Department of Mathematics and Statistics\\ Bar Ilan University, Ramat Gan, 52900, Israel\\        
}                                       



\dedication{Dedicated to the memory of J\"urgen Moser}            

\maketitle


 \begin{abstract}

 Multiplicative relations in the cohomology ring of a manifold impose
constraints upon its stable systoles.  Given a compact Riemannian
manifold $(X,g)$, its real homology $H_*(X, \RR)$ is naturally endowed
with the stable norm.  Briefly, if $h\in H_k (X, \RR)$ then the stable
norm of $h$ is the infimum of the Riemannian $k$-volumes of real
cycles representing $h$.  The stable $k$-systole is the minimum of the
stable norm over nonzero elements in the lattice of integral classes
in $H_k (X,\RR)$.  Relying on results from the geometry of numbers due
to W.~Banaszczyk, and extending work by M.~Gromov and J.~Hebda, we
prove metric-independent inequalities for products of stable systoles,
where the product can be as long as the real cup length of $X$.

 \end{abstract}






\section{Introduction}
Given a compact Riemannian manifold $(X,g)$, its real homology $H_*(X,
\RR)$ is naturally endowed with the stable norm.  Briefly, if $h\in
H_k (X, \RR)$ then the stable norm of $h$ is the infimum of the $k$-volumes,
defined in terms of the metric $g$, of real
cycles representing $h$.  The stable $k$-systole is the minimum of the
stable norm over nonzero elements in the lattice of integral classes
in $H_k (X,\RR)$, cf.~(\ref{(1.6)}) below.

Our goal is to explore the hidden power of the seminal calculation
contained in lines 7-12, p.~262 of M. Gromov's book \cite{[Gro3]},
cf.~(\ref{seminalgromov}) below.  We rely on results of W. Banaszczyk
\cite{[Ban1],[Ban2]} from the geometry of numbers, which provide a
bound for certain products of the successive minima of a pair of dual
lattices, see (\ref{(4.1)}) and (\ref{(7.1)}).  We thus extend the
work of M.~Gromov \cite{[Gro1]} and J.~Hebda \cite{[H]} to prove
metric-independent inequalities for products of stable systoles, where
the product can be as long as the real cup length of $X$.  In
particular, for the pair of stable systoles of dimension and
codimension 1, we prove an optimal inequality, where the boundary case
of equality is attained, e.g.\ by flat tori defined as quotients of the
so-called dual-critical lattices of A.-M.~Berg\'e and J. Martinet
\cite{[BM]}.

Geometric Measure Theory provides a framework in which the existence
of minimal representatives in homology classes can be proved.  If
$h\in H_k(X,\ZZ)$ then $\vol_k(h)$ (cf. Definition \ref{1.1} below) is
the minimal mass of a closed, integer-multiplicity, rectifiable
$k$-current representing $h$, cf.\ \cite[5.1.6]{[Fe1]}.  If $h\in
H_k(X,\RR)$ then the stable norm (or mass norm) of $h$ (cf.\
Definition \ref{1.2} below) is the minimal mass of a closed (normal)
$k$-current representing $h$.

In the cases $k=1$ and $k=n-1$ the transition from the Riemannian
length resp.~$(n-1)$-volume functional to the corresponding stable
norm is closely related to the process of {\it homogenization} studied
in Analysis.  For $k=n-1$ this shows up in J.~Moser's work \cite[Sect.~7]{[Mo]}, \cite{[MS]}, see \cite{[Se 1],[Se 2]} for more details.  It is through these
remarks by J.~Moser that the first-named author originally came in
contact with the subject of this paper.

While the present paper studies systolic ramifications of
multiplicative relations in the cohomology ring, the work \cite{[KKS]}
explores the systolic influence of Massey products.

Let $X$ be a compact manifold of dimension $n$.  We will now provide
detailed definitions of the systolic invariants involved.  A choice of
a Riemannian metric $g$ on $X$ allows us to define the total volume
$\vol_n(g)$, as well as the $k$-volumes of $k$-dimensional
submanifolds of $X$.  More generally, given an integer Lipschitz chain
$c=\sum_i n_i \sigma_i$, one defines its volume
$$
\vol_k(c)=\Sigma_i |n_i| \vol_k (\sigma_i).
$$
Here the volume $\vol_k(\sigma)$ of a Lipschitz $k$-simplex
$\sigma:\Delta^k\rightarrow X$ is the integral over the $k$-simplex
$\Delta^k$ of the ``volume form'' of the pullback $\sigma^*(g)$.

\begin{definition}\label{1.1}\rm
The (minimal) volume,
$$
\vol_k(h),
$$
of an integer homology class $h\in H_k(M,\ZZ)$ is the infimum of
$\vol_k(c)$ over all integer Lipschitz cycles representing $h$.
\end{definition}

Let $k$ be an integer satisfying $1\le k\le n$.  We define the
$k$-systole $\sys_k(g)$ of $(X,g)$ as the minimum of the volumes of
nonzero integer $k$-homology classes:
$$
\sys_k(g)=\min \{\vol_k(h)| h\in H_k(X, \ZZ)\setminus\{0\}\}.
$$

In particular, $\sys_n(g)=\vol_n(g)$ and, formally, $\sys_k(g)=\infty$
if $H_k(X,\ZZ)=0$.  A good introduction to systoles is M.~Berger's
survey \cite{[Be2]}, for further references see \cite{[KS1]}, sections 2 and 3.

Replacing integer cycles by real cycles we define the stable norm as
follows.

\begin{definition}\label{1.2}\rm
The stable norm $\parallelh$ of $h\in
H_k(X,\RR)$ is the infimum of the volumes $\vol_k(c)=\Sigma_i |r_i|
\vol_k(\sigma_i)$ over all real Lipschitz cycles $c=\Sigma_i r_i
\sigma_i$ representing $h$.
\end{definition}
Note that $\parallelleer$ is indeed a norm, cf.~\cite{[Fe2]} and \cite{[Gro3]},
4.C.

We denote by $H_k(X,\ZZ)_{\rr}$ the image of $H_k(X,\ZZ)$ in
$H_k(X,\RR)$ and by $h_{\rr}$ the image of $h\in H_k(X,\ZZ)$ in
$H_k(X,\RR)$.  Recall that $H_k(X,\ZZ)_{\rr}$ is a lattice in
$H_k(X,\RR)$.  Obviously
\begin{equation}
\label{(1.3)}
\parallelhrr \le \vol_k(h)
\end{equation}
for all $h\in H_k(X,\ZZ)$.  Moreover, $\parallelhrr = \vol_n(h)$ if
$h\in H_n(X,\ZZ)$.  H.~Federer \cite{[Fe2]} investiga\-ted the relations
between $\parallelhrr$ and $\vol_k(h)$ and proved:
\begin{equation}
\label{(1.4)}\quad\mbox{If $h\in H_k(X,\ZZ)$, $1\le k < n$, then}\hfill
\end{equation}
$$
\parallelhrr = \lim\limits_{i\rightarrow\infty} \frac{1}{i} \vol_k (i h).
$$
\begin{equation}
\label{(1.5)}\quad
\mbox{If $X$ is orientable and $h\in H_{n-1}(X,\ZZ)$,
then}\hfill
\end{equation}
$$
\parallelhrr =  \vol_{n-1} (h)
$$
(see also \cite{[Wh]}).  We define the stable $k$-systole $\stsys_k(g)$ of
$(X,g)$ by
\begin{equation}
\label{(1.6)}\stsys_k(g) = \min\left\{\parallelh | h\in H_k(X,
\ZZ)_{\rr} \setminus\{0\}\right\}.
\end{equation}
If the $k$'th Betti number $b_k(X) = \dim H_k(X,\RR)$ of $X$ is
positive and $H_k(X,\ZZ)$ is free abelian, then (\ref{(1.3)}) implies
$$
\stsys_k(g)\le \sys_k(g).
$$
If $X$ is orientable then
$$
\stsys_n(g) = \sys_n(g) = \vol_n(g)
$$
and (\ref{(1.5)}) implies
$$
\stsys_{n-1}(g) = \sys_{n-1}(g).
$$

\section{Statement of main results}
 Assuming that the fundamental class of a compact, oriented manifold
$X$ can be written as a cup product of classes of dimensions
$k_1,\ldots, k_m$, M.~Gromov proved an upper bound for the product of
the stable systoles $\Pi_j \stsys_{k_j}(g)$ in terms of the volume of
$X$, see \cite[7.4.C]{[Gro1]} and also \cite[4.38]{[Gro3]}.  We extend
this result, provide a simpler proof and analyze the dependence of the
constant on the Betti numbers involved.

\begin{theorem}\label{TheoremA}
Let $X$ be a compact manifold and let $k\ge 1$
be an integer such that $H^k(X,\RR)$ is not zero and spanned by cup
products of classes of dimensions $k_1, \ldots, k_m$, $\Sigma_j k_j=k$.
Then, for every Riemannian metric $g$ on $X$, we have
$$
\prod\limits^{m}_{j=1} \stsys_{k_j}(g) \le C(k)
\left(
\prod\nolimits_j b_{k_j}(X)(1+\log b_{k_j}(X))
\right)
\stsys_k(g)
$$
 for a constant $C(k)$ only depending on $k$.
\end{theorem}

Using Poincar\'{e} duality, J.~Hebda \cite[Proposition 6]{[H]} bounds
the product of the stable systoles in complementary dimensions by the
volume.  We generalize his result as follows:

\begin{theorem}\label{TheoremB}
 Let $X$ be a compact manifold and let $p, q$ be
integers, $p+q\le\dim X$.  Suppose there exists $h\in
H_{p+q}(X,\ZZ)_{\rr}$ such that the cap product with $h$ induces an
injective map
$$
\alpha\in H^p (X,\ZZ)_{\rr}\rightarrow h\cap\alpha\in H_q(X,\ZZ)_{\rr}.
$$
Then, for every Riemannian metric $g$ on $X$, we have
$$
\stsys_p(g) \stsys_q(g) \le C(p,q) b_p (X)\!\parallelh
$$
for a constant $C(p,q)$ only depending on $p$ and $q$.  If our
assumption is satisfied for all $h\in
H_{p+q}(X,\ZZ)_{\rr}\setminus\{0\}$ then
\begin{equation}
\label{(2.1)}\stsys_p(g) \stsys_q(g) \le C(p,q) b_p (X)
\stsys_{p+q}(g).
\end{equation}
\end{theorem}

Examples of flat tori show that an inequality of type (\ref{(2.1)})
has to depend on the Betti number $b_p$ linearly.  We do not know if
in Theorem \ref{TheoremA} the dependence on the Betti numbers is
optimal.  Explicit values for the constants $C(p,q)$ can be computed.
Unless $\{p,q\} = \{1, \dim X-1\}$ the sharp constants $C(p,q)$ are
unknown to us.  In the case $\{p,q\} = \{1, \dim X-1\}$ we have the
following sharp result (Corollary \ref{CorollaryC} below) that
generalizes \cite[Theorem A]{[H]} to the case $b_1(X)>1$.  Given a
lattice $L$ in euclidean space, we set
$$
\lambda_1(L) = \min\{|v|\mid v\in L \setminus \{0\}\}.
$$
Let $b\in\NN$, and consider the Berg\'{e}-Martinet constant
$\gamma'_b$, see \cite{[BM]},
$$
\gamma'_b = \sup\{\lambda_1(L) \lambda_1(L^*) \mid L \mbox{ a lattice
in euclidean space $\RR^b$}\}
$$
where $L^*$ is the lattice dual to $L$.  Thus, the constant
$\gamma'_b$ is bounded above by the Hermite constant $\gamma_b$,
cf.~e.g.~\cite{[LLS]}, p.~334, and satisfies $\gamma'_1=1$ and the inequalities
\begin{equation}
\label{(2.2)}\gamma'_b \le \gamma_b \le \frac{2}{3} \: b \quad
\mbox{for all $b\ge 2$}
\end{equation}
and
\begin{equation}
\label{(2.3)}
\frac{b}{2\pi e} (1+o (1)) \le \gamma'_b \le \frac{b}{\pi
e} \: (1+o (1)) \quad \mbox{for $b\rightarrow\infty$,}
\end{equation}
cf.\ e.g.~\cite{[LLS]}, pp.~334 and 337.

\begin{corollary}\label{CorollaryC}
Let $X$ be a compact, orientable manifold,
$\dim X=n$, with positive first Betti number $b=b_1(X)$.  Then, for
every Riemannian metric $g$ on $X$, we have
$$
\stsys_1(g) \sys_{n-1}(g) \le \gamma'_b \: \vol_n(g).
$$
 Equality is attained for a flat torus $\RR^n / L$ where
$L\subseteq \RR^n$ is a lattice with $\lambda_1(L) \lambda_1 (L^*) =
\gamma'_n$.
\end{corollary}

These inequalities can be seen as analogues of the optimal
inequalities of C. Loewner, P. Pu  \cite{[Pu]}, and C. Bavard \cite{[Bav]}, and as
generalizations of the results of R.~Accola \cite{[Ac]} and C.~Blatter
\cite{[Bl]} (cf.~\cite{[Be2]}, p.~290) in dimension 2.  A relative
version of Corollary \ref{CorollaryC} is studied in \cite{[Bab3]}.

Note that we have examples $X$ with equality in Corollary
\ref{CorollaryC} only if $\dim X =n\ge b=\dim H_1(X,\RR)$.  If we fix
$n=\dim X$ there might exist a better estimate, maybe even one
independent of $b$.  In particular we ask, see also \cite{[Gro1]},
7.4.C:

\begin{question}\label{2.4}\rm
Does there exist a constant $C$ such that
$$
\stsys_1(g) \sys_2(g) \le C\vol_3(g)
$$
for all 3-dimensional, compact, orientable, Riemannian manifolds $(X,g)$?
\end{question}

The manifolds $(X,g)$ for which equality holds in Corollary
\ref{CorollaryC} will be investigated in \cite{[BK]}.  In the case
$b_1(X)=1$, J.~Hebda \cite[Theorem~A]{[H]} proved that one has
equality $\stsys_1(g)\sys_{n-1}(g)=\vol_n(g)$ if and only if $X$
admits a Riemannian submersion $F:X\rightarrow S^1$ onto a circle,
such that all fibers $F^{-1}(s)$, $s\in S^1$, are minimal
hypersurfaces in $X$.

Consider a $b$-dimensional normed real vector space $(V,
\parallelleer)$ and a lattice $L\subseteq V$, in particular
$\rk(L)=b$.  For $1 \le i\le b$ we define the $i$-th successive
minimum $\lambda_i (L, \parallelleer)$ as the minimal $\lambda>0$ for
which there exist $i$ linearly independent vectors in $L$ of norm
smaller than or equal to $\lambda$. 

Now let $n=2p$, and consider the $L^2$-norm $\strichleer_{L^2}$ in
homology $H_p(X,\RR)$ dual to the one on harmonic forms defined by
$$
|f|_{L^2}^2=\int_X f\wedge *f,
$$
where $*$ is the Hodge star operator of the metric $g$.  Consider also
the conformally invariant norm $\|\;\|^*_{L^2}$ in $H^p_{dR}(X)$
defined by taking the infimum, over all representatives $\omega$, of
the quantity
$$
\parallelom_{L^2}= \left( \int_X \parallelomx^2 d\vol_n(x)
\right)^{\frac{1}{2}},
$$
where $\|\omega_x\|$ is the pointwise comass.  Let $\|\;\|_{L^2}$ be
the dual norm in homology, cf.\ (\ref{(7.4)}) below.  We introduce the
conformal invariant
$$
\conf_p(g)=\min\{\parallelh_{L^2} \;\mid h\in H_p
(X,\ZZ)_{\rr}\setminus\{0\}\}.
$$

\begin{corollary}\label{CorollaryD}
Let $p=\frac{n}{2}$ {\it and} $b_p=b_p(X)$.
Then
$$
\conf_p(g)^2 \le \frac{\lambda_1}{\lambda_{b_p}} {n\choose p}
b_p
$$
where $\lambda_i = \lambda_i(H_p(X,\ZZ)_{\rr}, \strichleer
_{L^2})$.
\end{corollary}

Note that 
\cite{[K2]} studies {\it lower} bounds for the quantity
$\sup_g\conf_2(g)$ for 4-manifolds $(X,g)$, where $X$ is the blowup of
the complex projective plane, the bound being {\it polynomial} of order
$n^{1\over 4}$ in the number $n$ of blown-up points.  Whether or not
such lower bounds exist for the stable 2-systole is unknown.

\begin{question}\label{2.5}\rm
Does there exist a constant $C$ such that
$$
\stsys_2(g)  \le C\sqrt{\vol_4(g)}
$$
for all 4-dimensional, compact, orientable, Riemannian manifolds $(X,g)$?
\end{question}

Note that this is so in the case of a definite intersection form,
cf.~\cite{[K2]}, (3.5).

For the first stable systole M.~Gromov proved
\begin{equation}
\label{gromsharp}
\stsys_1(g)\le\sqrt{\gamma_n}\vol_n(g)^{\frac{1}{n}}
\end{equation}
if $g$ is a Riemannian metric on a compact Riemannian manifold $X$
such that $b_1(X)=n=\dim X$, for which there are classes $\alpha_1,
\ldots, \alpha_n\in H^1(X,\RR)$ with
$\alpha_1\cup\ldots\cup\alpha_n\not= 0$.  Here $\gamma_n$ is the
Hermite constant, cf.\ (\ref{(2.2)}) above.  See
\cite[1.A$_1$]{[Gro2]}, \cite{[Gro3]}, pp.~259-260, \cite{[Be2]},
p.~283 and additionally \cite{[Bab1]}, Theorem 8.8.  Combining the
methods used to prove Theorem \ref{TheoremA} and Theorem
\ref{TheoremB}, we obtain a result similar to (\ref{gromsharp}) for
higher-dimensional stable systoles:

\begin{theorem}\label{TheoremE}
Let $X$ be a compact, orientable manifold, of dimension $n=kp$.
Suppose there are classes $\beta_1, \ldots, \beta_k\in H^p(X,\RR)$
such that $\beta_1\cup\ldots\cup\beta_k\not= 0$.  Then, for every
Riemannian metric $g$ on $X$, we have
$$
\stsys_p(g)\le C(n)\left( b_p^{\phantom{a}}(X) (1+\log b_p(X))
\right)^{\frac{k-1}{k}} \vol_n(g)^{\frac{1}{k}}
$$
for a constant $C(n)$ only depending on $n$.
\end{theorem}

It is a natural question whether similar inequalities exist if one
replaces the stable systoles by systoles.  The most prominent example
is Gromov's \cite[3.$C_1$]{[Gro2]} systolic inequality
\begin{equation}
\label{gromn}
\sys_1(g) \le c_n(\vol_n(g))^{\frac{1}{n}}
\end{equation}
which holds for all Riemannian metrics on compact $n$-dimensional
manifolds $X$ of cuplength $n$, i.e.~if for some field $F$ there are
classes $\alpha_1, \ldots, \alpha_n \in H^1(X,F)$ such that
$\alpha_1\cup\ldots\cup\alpha_n\not=0$.  See (\ref{grom2}) below for a
refined inequality of this type for surfaces.

For $k$-systoles with $k\ge 2$, however, systolic freedom tends to
prevail, i.e.~there are many examples where inequalities as in Theorem
\ref{TheoremA} and Theorem \ref{TheoremB} are not true if one replaces
stable systoles by systoles.  Based on examples by M.~Gromov
\cite[4.A$_5$]{[Gro2]}, see also section 6 of M. Berger's survey
article \cite{[Be2]}, this phenomenon has been studied by I.~Babenko,
A.~Suciu and the second author, cf.~\cite{[BaK]}, \cite{[BaKS]},
\cite[Appendix D]{[Gro3]}, \cite{[Bab2]}, \cite{[KS1],[KS2]} and
\cite{[K1]}.  Finally there is the striking example by M.~Freedman
\cite{[Fr]} of a sequence of Riemannian metrics $g_j$ on $S^1\times
S^2$ such that $\lim_j \vol_3(g_j)=0$ while length$(\gamma)$
area$(\Sigma)\ge 1$ (for all $j$) whenever $\gamma$ is a
noncontractible loop in $X$ and $\Sigma\subset X$ is a closed surface
which is not nullhomologous with coefficients in $\ZZ/2\ZZ$.

Some topological preliminaries are presented in section 3.  The
inequality of W. Banaszczyk, a crucial ingredient in our technique,
appears in section 4, together with the proof of Theorem
\ref{TheoremA}.  Theorem \ref{TheoremB} is proved in section 5,
Corollary \ref{CorollaryC} in section 6, Corollary \ref{CorollaryD} in
section 7, and Theorem \ref{TheoremE} in section~8.

\section{Some topological preliminaries}
In this section we collect
some facts from homology theory that will be used in the sequel.


For a compact manifold $X$ we consider the homology modules derived
from the singular (co-)chain complexes with $\ZZ$- or
$\RR$-coefficients.
We assume that the singular simplices are Lipschitz.
Since we will often use results that are actually proved for the
(continuous) singular or the $C^{\infty}$-singular (co-\nolinebreak[4])ho\-mology
theories we note that these are naturally isomorphic to the ones based
on Lipschitz simplices.  The bilinear pairing (Kronecker product)
between homology and cohomology will be denoted by [ , ], cup products
by $\cup$ and cap products by $\cap$.  By de Rham's Theorem, see
e.g.~\cite[5.45]{[Wa]}, the singular cohomology algebra $H^*(X,\RR)$ is
naturally isomorphic to the de Rham cohomology algebra $H^*_{dR}(X)$.
In de Rham theory the cup product is defined by the wedge product of
closed forms representing the cohomology classes.  The (nondegenerate)
bilinear pairing with real homology is given by integrating
representing forms over representing Lipschitz cycles, cf.~\cite[4.17]{[Wa]}.

\begin{remark}\label{3.1}\rm
The maps induced from the universal coefficient theorem are compatible
with Kronecker pairing, cup and cap products.
\end{remark}

More specifically, the map $\iota_*: H_*(X,\ZZ)\rightarrow
H_*(X,\RR)$, $h\rightarrow h_{\RR}$ is defined by inclusion on the
level of chains.  The map $\iota^*: H^*(X,\ZZ)\rightarrow H^*(X,\RR)
\simeq H^*_{dR}(X) $, $\alpha \rightarrow \alpha_{\RR}$ corresponds to
the natural extension of a cochain $c\in\Hom (C(X,\ZZ), \ZZ)$ to a
cochain in $\Hom(C(X,\RR),\RR)$.  According to the universal
coefficient theorem in homology, cf.~\cite[29.12]{[Gre]}, the image
$H_*(X,\ZZ)_\RR:=\Imag(\iota_*)$ is a lattice in $H_*(X,\RR)$.  Let
$H^*(X,\ZZ)_{\rr}:=\Imag(\iota^*)$.
\begin{lemma}
\label{(3.2)} The lattice $H^*(X,\ZZ)_{\rr}$ is dual to
$H_*(X,\ZZ)_{\rr}$ under $[\ ,\ ]$.
\end{lemma}

\begin{proof}
By naturality mentioned in Remark \ref{3.1} above, we have $[h_{\rr},
\alpha_{\rr}]=[h,\alpha]\in\ZZ$ for all $h\in H_*(X,\ZZ)$, $\alpha\in
H^*(X,\ZZ)$.  Conversely, suppose $\tilde{\alpha}\in H^*(X,\RR)$ and
$[h_{\rr}, \tilde{\alpha}]\in\ZZ$ for all $h\in H_*(X,\ZZ)$.  By the
surjectivity result \cite[23.9]{[Gre]} for a PID coefficient ring,
there exists $\alpha\in H^*(X,\ZZ)$ such that
$[h,\alpha]=[h_{\rr},\tilde{\alpha}]$ for all $h\in H_*(X,\ZZ)$.  The
nondegeneracy of [ , ] now implies that $\tilde{\alpha} = \alpha_{\rr}
\in H^*(X,\ZZ)_{\rr}$.
\end{proof}

Thus the set of de Rham classes $\alpha\in H^*_{dR}(X)$, that are represented
by closed forms whose integrals over integer cycles are integers,
coincides with $\Imag(\iota^*)$. These classes will henceforth be called
\emph{integer classes}.

If $X$ is a compact, oriented manifold, $\dim X=n$, then the
Poincar\'{e} duality isomorphism $H^p(X,\RR)\rightarrow
H_{n-p}(X,\RR)$ is defined by
$$
\PD_{\rr}:\alpha\in H^p(X,\RR)\rightarrow \xi_{\rr}\cap\alpha\in
H_{n-p}(X,\RR)
$$
where $\xi$ denotes the fundamental class of $X$.  Since the
corresponding map $\PD:\alpha\in H^p(X,\ZZ)\rightarrow\xi\cap\alpha\in
H_{n-p}(X,\ZZ)$ is an isomorphism, cf.~\cite[26.6]{[Gre]}, we see that
$\PD_{\rr}$ induces an isomorphism between $H^p(X,\ZZ)_{\rr}$ and
$H_{n-p}(X,\ZZ)_{\rr}$.  In particular, Lemma \ref{(3.2)} and the
preceding argument show that a de Rham class $\alpha\in H^p_{dR}(X)$ is
an integer class if $[\PD_{\rr}(\alpha),\beta]\in\ZZ$ for all integer
classes $\beta\in H^{n-p}_{dR}(X)$.  This fact was introduced as a
hypothesis in \cite{[Be1]}, p.~253, and named a ``dual lattice
condition'' in \cite{[H]}, p.~344.  The preceding discussion shows
that this condition is always satisfied.

\section{Banaszczyk's inequality and proof of Theorem~\ref{TheoremA}}
The following result by W.~Banaszczyk \cite{[Ban2]} is a crucial
ingredient in the proof of Theorem \ref{TheoremA}.

Consider a $b$-dimensional normed real vector
space $(V, \parallelleer)$ and a lattice $L\subseteq V$, in particular
$\rk(L)=b$.  For $1 \le i\le b$ we let $\lambda_i (L, \parallelleer)$
denote the minimal $\lambda>0$ for which there exist $i$ linearly
independent vectors in $L$ of norm smaller than or equal to $\lambda$.
Let $L^*$ denote the lattice dual to $L$ in the dual space $V^*$ and
let $\parallelleer^*$ denote the norm dual to $\parallelleer$.
Corollary 2 in \cite{[Ban2]} implies:

There exists a constant $C>0$ such that
\begin{equation}
\label{(4.1)} \lambda_i(L,\parallelleer)\lambda_{b-i+1}
(L^*,\parallelleer^*) \le C \, b(1+\log b)
\end{equation}
whenever $b\in\NN$, $i\in\{1, \ldots, b\}$ and $L$ is a lattice in a $b$-dimensional normed space
$(V, \parallelleer)$.

In our application of this result $V$ will be $H_p(X,\RR)$, $L$ will
be $H_p(X,\ZZ)_{\rr}$ and $\parallelleer$ will be the stable norm on
$H_p(X,\RR)$.  So $V^*$ can be identified with the de Rham cohomology
$H^p_{dR}(X)$, and then Lemma \ref{(3.2)} implies that
$L^*=H^p(X,\ZZ)_{\rr}$ is the lattice of integer de Rham classes.  The
norm on $H^p_{dR}(X)$ dual to the stable norm on $H_p(X,\RR)$ is the
comass norm, cf.~\cite[4.10]{[Fe2]} and \cite[4.35]{[Gro3]}.  Here the
comass of a $p$-form $\omega$ on a compact Riemannian manifold $X$ is
$$
\parallelom_{\infty}=\max\{\omega_x(e_1, \ldots, e_p) \mid x\in M,
e_i\in \TM_x \},
$$
where $|e_i|=1$ for $1\le i\le p$, while the comass of $\alpha\in
H^p_{dR}(X)$ is
$$
\parallelalf^* = \inf\{\parallelom_{\infty}\mid \omega \mbox{ a closed
$p$-form representing $\alpha$} \}.
$$
\begin{example}\label{4.2}\rm
 For a flat torus $(T^n, g)$ the calculation of the stable norm and
the comass norm reduces to the calculation of the corresponding
pointwise quantities, cf.~the work \cite{[Law]} by H.B.~Lawson, in
particular Corollary~4.4.  If we normalize $g$ so that $\vol_n(g)=1$,
and represent $(T^n, g)$ as the quotient of euclidean space $\RR^n$ by
a lattice $L$, then the standard isomorphisms
\begin{eqnarray*}
H_p(T^n, \RR)\simeq\Lambda^p \RR^n\\
H^p(T^n, \RR)\simeq\Lambda^p \RR^n
\end{eqnarray*}
convert the stable norm on $H_p(T^n,\RR)$ into the ``mass norm'' on
$\Lambda^p \RR^n$ (induced by the euclidean structure on $\RR^n$), and
the comass norm on $H^p(T^n, \RR)$ into the comass norm on $\Lambda^p
\RR^n$.  The cases $p=1$ and $p=n-1$ are particularly simple, since
then the mass and comass norms on $\Lambda^1 \RR^n\simeq \Lambda^{n-1}
\RR^n\simeq\RR^n$ coincide with the euclidean norm.  Moreover the
lattices $H_1(T^n,\ZZ)_{\rr}\subset H_1(T^n,\RR)\simeq\RR^n$ and
$H_{n-1}(T^n,\ZZ)_{\rr} \subset H_{n-1}(T^n,\RR)\simeq\RR^n$
correspond to $L$ resp.~$L^*$ under these isomorphisms.  The cases
$1<p<n-1$ are considerably more complicated, see \cite{[Law]}.
\end{example}

We are now in a position to prove our first theorem.

\begin{otheorem}{\bf Theorem \ref{TheoremA}}
Let $X$ be a compact manifold and let $k\ge 1$ be an integer such that
$H^k(X,\RR)$ is not zero and spanned by cup products of classes of
dimensions $k_1, \ldots, k_m$, $\Sigma_j k_j=k$.  Then, for every
Riemannian metric $g$ on $X$, we have
$$
\prod\limits^{m}_{j=1} \stsys_{k_j}(g) \le C(k)
\left(
\prod\nolimits_j b_{k_j}(X)(1+\log b_{k_j}(X))
\right)
\stsys_k(g)
$$
 for a constant $C(k)$ only depending on $k$.
\end{otheorem}

\begin{proof}[Proof of Theorem \ref{TheoremA}]
For notational reasons we present the proof for the case that
$H^k(X,\RR)$ is non-zero, and is spanned by cup products of classes of
dimensions $j$ and $l$, with $j+l=k$.  The generalization to more
factors is obvious.  Since $H^k(X,\RR)\not=\{0\}$, there exists $h\in
H_k(X,\ZZ)_{\rr}$ such that
\begin{equation}
\label{(4.3)}\parallelh=\stsys_k(g).
\end{equation}
We set $\lambda^* = \lambda_{b_j}(H^j(X,\ZZ)_{\rr}, \parallelleer^*)$
and $\mu^* = \lambda_{b_l}(H^l(X,\ZZ)_{\rr}, \parallelleer^*)$.  Then
we can find integral classes $\alpha_1, \ldots, \alpha_{b_j}$ in
$H^j(X,\ZZ)_{\rr}$ spanning $H^j(X,\RR)$ and $\beta_1, \ldots,
\beta_{b_l}$ in $H^l(X,\ZZ)_{\rr}$ spanning $H^l(X,\RR)$ such that
\begin{equation}
\label{(4.4)}\mbox{
$\parallelalfs^* \le \lambda^*$
for
$1\le s\le b_j$,
and
$\parallelbetat^* \le\mu^*$
for
$1\le t\le b_l$.
}
\end{equation}
By assumption $H^k(X,\RR)$ is spanned by the cup products $\alpha_s
\cup \beta_t$, $1\le s\le b_j$ and $1\le t\le b_l$.  Hence there exist
indices $s$ and $t$ such that
$$
[h, \alpha_s\cup \beta_t] = a\not= 0.
$$
According to Remark \ref{3.1}, we have $a\in\ZZ$.  Now we complete the
proof by a calculation analogous to Gromov's calculation in
\cite[Theorem~4.36]{[Gro3]}, p.~262. If $c$ is a real Lipschitz cycle
representing $h$ and if $\omega$ and $\pi$ are closed forms
representing $\alpha_s$ and $\beta_t$, then
\begin{equation}
\label{seminalgromov}
1\le |a| = \left|\int_c \omega\wedge\pi\right| \le \frac{k!}{j!\, l!}
\vol_k(c) \parallelom_{\infty} \parallelpi_{\infty},
\end{equation}
cf.~\cite[1.8.1]{[Fe1]} for the factor $\frac{k!}{j!\, l!} = {k\choose j}$.
Using this and (\ref{(4.3)}) and (\ref{(4.4)}) we obtain
$$
1\le |a| \le \frac{k!}{j!\, l!} \stsys_k(g) \lambda^* \mu^*.
$$
Since $H^p(X,\ZZ)_{\rr}$ is the lattice dual to $H_p(X,\ZZ)_{\rr}$ by
Lemma \ref{(3.2)}, we can use Banaszczyk's inequality (\ref{(4.1)}) to
conclude
$$
\stsys_j(g) \stsys_l(g) \le \frac{k!}{j!\, l!} C^2 \,
b_j(1+\log b_j)
b_l(1+\log b_l)
\stsys_k(g).
$$
This is the claim of Theorem \ref{TheoremA} in the case $m=2$.  For
arbitrary $m\ge 2$ we obtain the constant
$$
\frac{k!}{\Pi_j \, k_j!} C^m
$$
in this inequality.
\end{proof}

\section{Inequalities based on cap products and Poincar\'{e} duality}
In this section we prove Theorem \ref{TheoremB}.  The proof of our
Theorem~\ref{TheoremA} was based on Banaszczyk's inequality
(\ref{(4.1)}) applied to the stable norm.  For Euclidean norms,
Banaszczyk proved a sharper estimate which is linear in the dimension,
cf.\ (\ref{(7.1)}) below.  Much older, linear estimates for the
Hermite constant itself were already mentioned above, cf.\
(\ref{(2.2)}) and (\ref{(2.3)}).  Replacing the stable norm by the
$L^2$-norm, J.~Hebda \cite{[H]} was able to apply such estimates in
complementary dimensions: if $\vol_n(g)$ is normalized to one, then
the $L^2$-norm on forms is bounded above by (a constant times) the
comass norm and consequently the stable norm is bounded above by (a
constant times) the dual $L^2$-norm.  A simple estimate shows that the
stable norm of the Poincar\'{e} dual of a cohomology class is bounded
above by (a constant times) its $L^2$-norm.  Putting these facts
together we see that every inequality of type (\ref{(4.1)}) leads to a
systolic inequality in complementary dimensions.  In Theorem
\ref{TheoremB} we generalize this procedure, by replacing the
fundamental cycle used in Poincar\'{e} duality by a cycle of arbitrary
dimension.

First we discuss a generalization of the $L^2$-norm that is
appropriate in this context.  We stay in the well-known realm of
Lipschitz cycles, although at some points the presentation would be
more elegant if we used closed currents of finite mass.

If $z=\Sigma_i\, r_i \sigma_i$ is a real Lipschitz $(p+q)$-cycle in
the compact Riemannian manifold $(X,g)$, we define a positive
semidefinite, symmetric bilinear form $\langle$ , $\rangle_{L^2 (z)}$
on the space of $p$-forms as follows.  Let $g_i= \sigma_i^* g$ denote
the pull-back of $g$ by the Lipschitz simplex $\sigma_i$. So $g_i$ is
a positive semidefinite, symmetric (0,2)-tensor field that is bounded,
measurable and defined almost everywhere on the standard
$(p+q)$-simplex $\Delta^{p+q}$.  Let $\langle,\rangle_{g_i}$ denote
the bilinear form induced by $g_i$ on the bundle of alternating
$p$-tensors over $\Delta^{p+q}$.  Let $d\vol_{g_i}$ denote the
``volume element'' induced by $g_i$ on $\Delta^{p+q}$, in particular
$\vol_{g_i}(\Delta^{p+q})=\vol_{p+q}(\sigma_i)$.  If $\omega$ and
$\overline{\omega}$ are $p$-forms on $M$ we set
$$
\langle\omega,\overline{\omega}\rangle_{L^2(z)} = \sum_i |r_i|
\int_{\Delta^{p+q}}
\langle\sigma_i^*\omega, \sigma_i^*\overline{\omega}\rangle_{g_i}
d\vol_{g_i}.
$$
Since $\langle(\sigma_i^*\omega)_x, (\sigma_i^*\omega)_x \rangle_{g_i}
\le {p+q \choose p} {\norm {\omega_{\sigma_i (x)}}}^2$ for
$x\in\Delta^{p+q}$, cf.\ \cite[1.8.1]{[Fe1]}, we have
\begin{equation}
\label{(5.1)}\langle\omega,\omega\rangle_{L^2(z)} \le {p+q \choose p}
\parallelom^2_{\infty} \vol_{p+q}(z).
\end{equation}
More generally, if $\pi$ is an additional $q$-form, we can estimate
\begin{equation}
\label{(5.2)}\left(\int_z \omega\wedge\pi\right)^2 \le \, {p+q \choose
p} \parallelom^2_{\infty} \langle\pi,\pi\rangle_{L^2(z)}
\vol_{p+q}(z).
\end{equation}
Obviously, we can interchange the roles of $\omega$ and $\pi$ in this inequality.
If $p+q<n=\dim X$, there will be many $p$-forms $\omega$ such that
the support of
$\omega$ is disjoint from all $\sigma_i(\Delta^{p+q})$, and
hence $\langle\omega,\omega\rangle_{L^2(z)}=0$.  However, under
appropriate homological conditions on $z$ the semidefinite form
$\langle\ ,\ \rangle_{L^2(z)}$ induces a scalar product on
$H^p(X,\RR)$.

\begin{lemma}\label{5.3}
Suppose $h\in H_{p+q}(X,\RR)$ satisfies $h\cap
\alpha\not= 0$ for all $\alpha\in H^p(X,\RR)\setminus\{0\}$.  If $z$
is a Lipschitz cycle representing $h$, then there exists a scalar
product on $H^p(X,\RR)$, denoted by
$\langle$ , $\rangle_{L^2(z)}$ as well, such that for all $\alpha\in
H^p(X,\RR)$
$$
\langle\alpha,\alpha\rangle_{L^2(z)} = \inf \left\{
\langle\omega,\omega\rangle_{L^2(z)} \mid \omega \mbox{ a closed
$p$-form representing $\alpha$}\right\}.
$$
\end{lemma}
\begin{proof}
The form $\langle$ , $\rangle_{L^2(z)}$ on the space $Z^p(X)$
of closed $p$-forms descends to a scalar product $\langle$ ,
$\rangle_{L^2(z)}$ on $Z^p(X)/N$, where
$$
N=\left\{\omega\in Z^p(X) \mid \langle\omega,\omega\rangle_{L^2(z)}
=0\right\}.
$$
Now $N$ is contained in the space $B^p(X)$ of exact forms. Indeed, if
$\omega\in Z^p(X)\setminus B^p(X)$ represents $\alpha\in
H^p(X,\RR)\setminus\{0\}$ then $h\cap\alpha\not= 0$ by assumption.
Hence there exists $\beta\in H^p(X,\RR)$ such that
$[h\cap\alpha,\beta]\not= 0$.  It follows from \cite[24.19]{[Gre]} that
$[h\cap\alpha,\beta]=[h,\alpha\cup\beta]$.  If $\pi$ is a closed
$q$-form representing $\beta$, then by (\ref{(5.2)})
$$
0<[h,\alpha\cup\beta]^2 \le {p+q \choose p} \parallelpi^2_{\infty}
\langle\omega,\omega\rangle_{L^2(z)} \vol_{p+q}(z)
$$
and hence $\langle\omega,\omega\rangle_{L^2(z)} >0$.

Arguing in the completion of $(Z^p(X)/N$, $\langle$ ,
$\rangle_{L^2(z)})$ and using the fact that
$Z^p(X)/B^p(X)\simeq(Z^p(X)/N)/(B^p(X)/N)$ is finite-dimensional, one
can easily complete the proof.
\end{proof}

We are now in a position to prove our second theorem.

\begin{otheorem}{\bf Theorem \ref{TheoremB}}
Let $X$ be a compact manifold and let $p, q$ be integers, $p+q\le\dim
X$.  Suppose there exists $h\in H_{p+q}(X,\ZZ)_{\rr}$ such that the
cap product with $h$ induces an injective map
$$
\alpha\in H^p (X,\ZZ)_{\rr}\rightarrow h\cap\alpha\in H_q(X,\ZZ)_{\rr}.
$$
Then, for every Riemannian metric $g$ on $X$, we have
$$
\stsys_p(g) \stsys_q(g) \le C(p,q) b_p (X)\!\parallelh
$$
for a constant $C(p,q)$ only depending on $p$ and $q$.  If our
assumption is satisfied for all $h\in
H_{p+q}(X,\ZZ)_{\rr}\setminus\{0\}$ then
$$
\tag{\mbox{\rm (\ref{(2.1)})}}{\stsys_p(g) \stsys_q(g) \le C(p,q) b_p (X)
\stsys_{p+q}(g).}
$$
\end{otheorem}

\begin{proof}[Proof of Theorem \ref{TheoremB}]
Let $z$ be a Lipschitz cycle representing $h$.  By the preceding lemma
we can consider the scalar product $\langle\; , \;\rangle_{L^2(z)}$ on
$H^p_{dR}(X)$.  From (\ref{(5.1)}) we obtain
\begin{equation}
\label{(5.4)}|\alpha|_{L^2(z)} \le {p+q \choose p}^{\frac{1}{2}}
\vol_{p+q}(z)^{\frac{1}{2}} \parallelalf^*
\end{equation}
for $\alpha\in H^p_{dR}(X)$.  Hence the dual euclidean norm on
$H_p(X,\RR)$, also denoted by $\strichleer_{L^2(z)}$, satisfies
\begin{equation}
\label{(5.5)}\parallelk \le {p+q \choose p}^{\frac{1}{2}}
\vol_{p+q}(z)^{\frac{1}{2}} |k|_{L^2(z)}
\end{equation}
for $k\in H_p(X,\RR)$.  If $\alpha\in H^p_{dR}(X)$ then
$$
\parallelhcapalf=\sup\left\{[h,\alpha\cup\beta] \mid \beta\in
H^q(X,\RR), \parallelbeta^*\le 1\right\}
$$
since the stable norm $\parallelleer$ and the comass norm $\parallelleer^*$
are dual to each other and $[h\cap\alpha,\beta]=[h,\alpha\cup\beta]$.
Now (\ref{(5.2)}) implies
\begin{equation}
\label{(5.6)}\parallelhcapalf \le {p+q \choose p}^{\frac{1}{2}}
\vol_{p+q}(z)^{\frac{1}{2}} \, |\alpha|_{L^2(z)}.
\end{equation}
Set $b=b_p(X)$.  By the definition of the constant $\gamma'_b$ and Lemma
\ref{(3.2)} there exist $k\in H_p(X,\ZZ)_{\rr}\setminus\{0\}$ and $\alpha\in
H^p(X,\ZZ)_{\rr}\setminus\{0\}$ such that
$$
|k|_{L^2(z)} \, |\alpha|_{L^2(z)} \le \gamma'_b.
$$
Using (\ref{(5.5)}) and (\ref{(5.6)}) we conclude
$$
\parallelk \, \parallelhcapalf \le {p+q \choose p} \gamma'_b \,
\vol_{p+q}(z).
$$
Since $h\cap\alpha\not= 0$ by assumption and $h\cap\alpha\in
H_q(X,\ZZ)_{\rr}$ by Remark \ref{3.1} we conclude
\begin{equation}
\label{(5.7)}\stsys_p(g) \, \stsys_q(g) \le {p+q \choose p} \gamma'_b
\parallelh.
\end{equation}
Together with (\ref{(2.2)}) this proves Theorem \ref{TheoremB}.
\end{proof}

\begin{remark}\label{5.8}\rm
In \cite{[H]}, Proposition 6, J.~Hebda essentially proved Theorem
\ref{TheoremB} in the case $n=p+q$. Indeed, his arguments imply that
$$
\stsys_p(g) \, \stsys_{n-p}(g)\le{n \choose p} \gamma_b \, \vol_n(g)
$$
for every Riemannian metric $g$ on a compact, orientable manifold $X$
with $b=b_p(X)>0$.  Theorem \ref{TheoremB}, however, can be applied in
many situations that are not covered by \cite[Proposition 6]{[H]}. In
general, the resulting inequality will be sharper than the inequality
provided by Theorem \ref{TheoremA}.  We give a simple example. We
choose integers $p\ge 1$, $q\ge 1$ and $m$ with $m>p+q$ and apply
(\ref{(5.7)}) to the manifold $X=T^{p+q}\times S^m$.  Recalling that
$\gamma'_b \le \frac{2}{3}b$ for $b\ge 2$, cf.\ (\ref{(2.2)}), we
obtain
\begin{equation}
\label{(5.9)}\stsys_p(g) \stsys_q(g) \le \frac{2}{3} {p+q \choose p}^2
\stsys_{p+q}(g)
\end{equation}
for every metric $g$ on $X=T^{p+q}\times S^m$.  Here, the dependence
on $b_p(X)=b_q(X)={p+q \choose p}$ is quadratic, while the inequality
following from the proof of Theorem \ref{TheoremA} is
$$
\stsys_p(g) \stsys_q(g) \le C^2 {p+q \choose p} b_p(X)^2
(1+\log b_p(X))^2 \stsys_{p+q}(g)
$$
where $C$ is the universal constant in Banaszczyk's inequality
(\ref{(4.1)}).
\end{remark}

\section{A sharp inequality in codimension 1}
We apply the arguments of section 5 to a fundamental cycle
$h=\xi_{\rr}\in H_n(X,\ZZ)_{\rr}$ of $X$.

\begin{otheorem}{\bf Corollary \ref{CorollaryC}}
Let $X$ be a compact, orientable manifold,
$\dim X=n$, with positive first Betti number $b=b_1(X)$.  Then, for
every Riemannian metric $g$ on $X$, we have
$$
\stsys_1(g) \sys_{n-1}(g) \le \gamma'_b \: \vol_n(g).
$$
 Equality is attained for a flat torus $\RR^n / L$ where
$L\subseteq \RR^n$ is a lattice with $\lambda_1(L) \lambda_1 (L^*) =
\gamma'_n$.
\end{otheorem}

\begin{proof}[Proof of Corollary \ref{CorollaryC}]
Since $p=1, q=n-1$ we can replace the factors ${p+q \choose
p}^{\frac{1}{2}}$ in (\ref{(5.5)}) and (\ref{(5.6)}) by one,
cf.~\cite[1.8.1]{[Fe1]}.  So our final statement is
\begin{equation}
\label{(6.1)}\stsys_1(g) \sys_{n-1}(g) \le \gamma'_b \, \vol_n(g).
\end{equation}
There exist lattices $L$ in $b$-dimensional euclidean space such that
$$
\lambda_1(L) \lambda_1(L^*)=\gamma'_b,
$$
cf.~\cite{[BM]}.  Consider the corresponding flat torus $(T^b,g)$.
Then Example~\ref{4.2} shows that for such tori equality holds in
(\ref{(6.1)}).
\end{proof}

The boundary case of equality in this optimal inequality will be
studied in \cite{[BK]}.

\section{A conformally invariant inequality in the middle dimension}
Instead of using the constant $\gamma'_b$ together with inequality
(\ref{(2.2)}) in the proof of Theorem \ref{TheoremB}, we can also
employ the following estimate by W.~Banaszczyk \cite{[Ban1]}, Theorem
2.1, in the case $i=b$:

If $L$ is a lattice in $b$-dimensional euclidean space (with the standard
innner product norm) and $i\in\{1, \ldots, b\}$, then
\begin{equation}
\label{(7.1)}1\le\lambda_i(L) \lambda_{b-i+1}(L^*)\le b.
\end{equation}
Specialized to the case $p+q=n$, $h=\xi_{\rr}$, the proof of Theorem
\ref{TheoremB}, combined with (\ref{(7.1)}), leads to the estimate
\begin{equation}
\label{(7.2)}\stsys_p(g) \, \stsys_{n-p}(g)\le
\frac{\lambda_1}{\lambda_{b_p}} {n\choose p} b_p \vol_n(g)
\end{equation}
where $\lambda_i=\lambda_i(H_p(X,\ZZ)_{\rr}, \strichleer_{L^2})$
and $b_p = b_p(X)$.

In the case of middle dimension $p=\frac{n}{2}$ a similar argument
leads to a conformally invariant estimate, cf.~\cite[7.4.A]{[Gro1]}.  This
case is special in the following sense: the $L^2$-norms of $p$-forms
$\omega$
$$
|w|_{L^2}=\left(\int_X \langle w,w\rangle d\vol_n
\right)^{\frac{1}{2}}
$$
and
$$
\parallelom_{L^2}= \left( \int_X ( \parallelomx)^2 d\vol_n(x)
\right)^{\frac{1}{2}}
$$
are invariant under conformal changes of the metric $g$.  They induce
respectively the usual $L^2$-norm $\strichleer_{L^2}$ on $H^p (X,\RR)$,
for which the harmonic forms are minimizing representatives of the
cohomology classes, and a norm $\parallelleer^*_{L^2}$ on $H^p (X,\RR)$
such that
$$
{n\choose p}^{-\frac{1}{2}} |\alpha|_{L^2} \le \, \parallelalf^*_{L^2}
\, \le \, |\alpha|_{L^2}.
$$
We also consider the conformally invariant dual norms on $H_p(X,\RR)$
that satisfy
\begin{equation}
\label{(7.3)} |h|_{L^2} \le \parallelh_{L^2} \le {n\choose
p}^{\frac{1}{2}} |h|_{L^2}.
\end{equation}
It is not difficult to prove that
\begin{equation}
\label{(7.4)}\parallelh_{L^2} = \sup
\{\parallelh_{g'} \vol_n(g')^{-\frac{1}{2}} \mid g' \mbox{
conformal to $g$}\},
\end{equation}
cf.~\cite[7.4.A]{[Gro1]}.
Always assuming $p=\frac{n}{2}$, we introduce the conformal invariant
$$
\conf_p(g)=\min\{\parallelh_{L^2} \mid h\in H_p
(X,\ZZ)_{\rr}\setminus\{0\}\}.
$$
Note that (\ref{(7.4)}) implies that
\begin{equation}
\label{(7.5)}
\stsys_p (g') \vol_n(g')^{-\frac{1}{2}} \le \conf_p(g)
\end{equation}
for every pair $g'$, $g$ of conformal Riemannian metrics.
However, it can happen that
$$
\sup\{\stsys_p(g') \vol_n(g')^{-\frac{1}{2}} \mid g' \mbox{ conformal
to $g$}\} < \conf_p(g).
$$
This can be seen by comparing Gromov's universal upper bound contained
in inequality (+) of \cite[2.$C$]{[Gro2]} for the 1-systole to the
lower bound for $\conf_1$ of surfaces which follows from the work by
P. Buser and P. Sarnak \cite{[BS]}.  Indeed, Gromov asserts the
existence of a universal constant $C$ such that
\begin{equation}
\label{grom2}
\sys_1(g) \vol_2(g)^{-\frac{1}{2}} \le C \,
\frac{\log\gamma}{\sqrt{\gamma}}
\end{equation}
whenever $g$ is a Riemannian metric on a closed, orientable surface of
genus $\gamma\ge 2$, cf.~(\ref{gromn}) above.  On the other hand
\cite[1.13]{[BS]} states that the supremum of $\conf_1(g)^2$ over all
Riemannian metrics on a closed, orientable surface of genus $\gamma$
is bounded below by $c\log(\gamma)$ for some universal constant
$c>0$. Thus, a priori, an upper bound for $\conf_p(g)$ as given in
Corollary \ref{CorollaryD} is stronger than the same upper bound for
$\stsys_p(g) \vol_n(g)^{-\frac{1}{2}}$, cf.~also \cite{[K2]}.

\begin{otheorem}{\bf Corollary \ref{CorollaryD}}
Let $p=\frac{n}{2}$ {\it and} $b_p=b_p(X)$.
Then
$$
\conf_p(g)^2 \le \frac{\lambda_1}{\lambda_{b_p}} {n\choose p}
b_p
$$
where $\lambda_i = \lambda_i(H_p(X,\ZZ)_{\rr}, \strichleer
_{L^2})$.
\end{otheorem}

\begin{proof}[Proof of Corollary \ref{CorollaryD}]
Note that the map
$$
\PD_{\rr} : (H^p(X,\RR), \strichleer_{L^2}) \rightarrow (H_{n-p}
(X,\RR), \strichleer_{L^2})
$$
is an iso\-me\-try. Indeed, if $\alpha\in H^p(X,\RR)$ and if $\omega$ is the
harmonic $p$-form representing $\alpha$, then
$$
\begin{array}{rcl}
|\PD_{\rr}(\alpha)|_{L^2} &=& \sup\{ \int_M \omega\wedge\pi\mid\pi\in
Z^{n-p} (M), |\pi|_{L^2}=1 \} \\ &=& \sup\{ \langle \omega, \,
*\pi\rangle_{L^2}\mid \pi\in Z^{n-p} (M), |\pi|_{L^2}= 1 \}
=|\alpha|_{L^2}.
\end{array}
$$
Setting $\lambda^*_1 = \lambda_1(H^p(X,\ZZ)_{\rr}, \strichleer_{L^2})$
and $b=b_p(X)$, we conclude $\lambda_1=\lambda^*_1$ and
$\lambda^2_1\le\frac{\lambda_1}{\lambda_b} b$ by Banaszczyk's
inequality (\ref{(7.1)}).  Now (\ref{(7.3)}) implies our claim.
\end{proof}

One can also apply inequality (\ref{(7.2)}) in the opposite direction:

\begin{corollary}\label{7.5}
Let $X$ be a compact, orientable manifold
of dimension $n=2p$, with $b_p(X)>0$. For $D>0$ let ${\mathcal G}_D$
denote the set of Riemannian metrics $g$ on $X$ such that
$\vol_{n}(g)\le D$ and $\stsys_p(g)\ge D^{-1}$.  Then the set of flat
Finsler metrics defined by the stable norms of metrics $g\in {\mathcal
G}_D$ on the Jacobian torus $H_p(X,\RR)/H_p(X,\ZZ)_{\rr}$ is contained
in a compact part of the set of all flat Finsler metrics on the
Jacobian torus.
\end{corollary}

\begin{proof}
First we will prove the existence of a basis of
$H_p(X,\ZZ)_{\rr}$ whose elements have stable norm bounded by some
function of $D$, $p$ and $b_p(X)$.  We use Banaszczyk's inequality
(\ref{(4.1)}) in the case $i=b=b_p(X)$ for the stable norm $\parallelleer$ on
$L=H_p(X,\ZZ)_{\rr}$ and the comass norm $\parallelleer^*$ on
$L^*=H^p(X,\ZZ)_{\rr}$ to obtain
\begin{equation}
\label{(7.6)}\lambda_b \lambda^*_1 \le C b(1+\log b)
\end{equation}
where $\lambda_i=\lambda_i(H_p(X,\ZZ)_{\rr}, \parallelleer),$
and $\lambda^*_i=\lambda_i(H^p(X,\ZZ)_{\rr}, \parallelleer^*)$.
If
$\alpha\in H^p (X,\ZZ)_{\rr}$ and $\parallelalf^*=\lambda^*_1$ then
(\ref{(5.4)}) and (\ref{(5.6)}) imply
$$
\parallelPDralf \le {n \choose p} \vol_{n}(g)\parallelalf^*
$$
and hence
$$
\stsys_p(g)\le {n \choose p} \vol_{n}(g) \lambda^*_1.
$$
Since $\stsys_p(g)=\lambda_1$, we can use (\ref{(7.6)}) to conclude
\begin{equation}
\label{(7.7)}\stsys_p(g)^2 \le \frac{\lambda_1}{\lambda_b} {n \choose
p} C \, b(1+\log b) \vol_{n}(g).
\end{equation}
Now assume $g\in \mathcal G_D$.  Then (\ref{(7.7)}) implies $\lambda_b \le {n
\choose p} C \, b(1+\log b)D^2$.  Using \cite{[C]}, p.135, Lemma 8, we obtain
a basis $v_1, \ldots, v_b$ of $H_p(X,\ZZ)_{\rr}$ such that
\begin{equation}
\label{(7.8)}\parallel v_i\parallel \le {n \choose p} C \, b^2(1+\log
b)D^2=E
\end{equation}
for $1\le i\le b$.  Consider the isomorphism
$H_p(X,\RR)\rightarrow\RR^b$ mapping $(v_1, \ldots, v_b)$ to the
standard basis $(e_1, \ldots, e_b)$ and consider the induced norm on
$\RR^b$, denoted by $\parallelleer$ as well.  Then the unit
$\parallelleer$-ball $B=\{x\in\RR^b \mid \parallelx\le 1\}$ satisfies
\begin{equation}
\label{(7.9)}\Big\{x\in\RR^b\mid \sum\limits^b_{i=1} |x_i| \le
E^{-1}\Big\} \subseteq B.
\end{equation}
On the other hand $\min\left\{\parallelx \mid x\in\ZZ^b\setminus
\{0\}\right\}=\lambda_1\ge D^{-1}$ and hence Min\-kows\-ki's Theorem,
cf.~\cite{[C]}, p.71, Theorem II, implies
\begin{equation}
\label{(7.10)}\vol_b(B)\le(2D)^b
\end{equation}
where $\vol_b$ denotes the usual Lebesgue measure on $\RR^b$.  Now our
claim follows from the fact that for fixed numbers $D$ and $E>0$ the
set of convex bodies in $\RR^b$ satisfying (\ref{(7.9)}) and
(\ref{(7.10)}) is compact with respect to the Hausdorff metric.
\end{proof}

\section{A sublinear estimate for a single systole}
If $X$ is a compact, orientable manifold, if $\dim X=n=kp$ and if
there are classes $\alpha_1, \ldots, \alpha_k \in H^p(X,\RR)$ such
that $\alpha_1\cup\ldots\cup\alpha_k\not= 0$, then Theorem
\ref{TheoremA} implies
$$
\stsys_p(g)\le\tilde{C}(n) \, b_p(X) (1+\log b_p(X))
\vol_n(g)^{\frac{1}{k}}
$$
for every Riemannian metric $g$ on $X$.  Using Poincar\'{e} duality as
in Theorem \ref{TheoremB} we can improve this estimate to sublinear
dependence on the Betti number $b=b_p(X)$.

\begin{otheorem}{\bf Theorem \ref{TheoremE}}
Let $X$ be a compact, orientable manifold, of dimension $n=kp$.
Suppose there are classes $\beta_1, \ldots, \beta_k\in H^p(X,\RR)$
such that $\beta_1\cup\ldots\cup\beta_k\not= 0$.  Then, for every
Riemannian metric $g$ on $X$, we have
$$
\stsys_p(g)\le C(n)\left( b_p^{\phantom{a}}(X) (1+\log b_p(X))
\right)^{\frac{k-1}{k}} \vol_n(g)^{\frac{1}{k}}
$$
for a constant $C(n)$ only depending on $n$.
\end{otheorem}
\begin{proof}[Proof of Theorem \ref{TheoremE}]
We start by choosing linearly independent classes $\alpha_1, \ldots,
\alpha_b$ in $H^p(X,\ZZ)_{\rr}$ such that $\parallelalfi^* \le
\lambda^*_b= \lambda_b(H^p(X,\ZZ)_{\rr}, \parallelleer^*)$.  Our
hypothesis implies that there exist indices $i_1, \ldots, i_k$ in
$\{1, \ldots, b\}$ such that
$$
\alpha_{i_1}\cup\ldots\cup\alpha_{i_k}\not= 0.
$$
We may assume that $X$ is connected and oriented. Let $h\in
H_p(X,\ZZ)_{\rr}\setminus\{0\}$ denote the Poincar\'{e} dual of
$\alpha_{i_1}\cup\ldots\cup\alpha_{i_{k-1}}$.  Let
$\omega_1,\ldots,\omega_b$ be $p$-forms representing
$\alpha_1,\ldots,\alpha_b$, and let $\omega$ be an arbitrary closed
$p$-form representing some $\alpha\in H^p(X,\RR)$.  Then
$$
\begin{array}{rcl}
[h,\alpha]&=&\left[\xi_{\rr}, \alpha_{i_1}\cup\ldots\cup
\alpha_{i_{k-1}}\cup\alpha\right] = \int_X
\omega_{i_1}\wedge\ldots\wedge \omega_{i_{k-1}}\wedge\omega\\ &\le&
\frac{n!}{(p!)^k} \left(\prod\limits^{k-1}_{j=1}
\parallelomij_{\infty}\right) \parallelom_{\infty} \vol_n(X).
\end{array}
$$
Since $\parallelh=\sup\limits_{\parallelalf^*\le 1}[h,\alpha]$, we obtain
$$
\parallelh\le \frac{n!}{(p!)^k} \left(\prod\limits^{k-1}_{j=1}
\parallelalfij^*\right) \vol_n(X)\le \frac{n!}{(p!)^k}
\left(\lambda^*_b\right)^{k-1} \vol_n(X).
$$
Now Banaszczyk's inequality (\ref{(4.1)}) implies $\stsys_p(g)
\lambda^*_b \le C b(1+\log b)$, and hence
$$
\stsys_p(g)^k \le \stsys_p(g)^{k-1} \parallelh \le C^{k-1}
\frac{n!}{(p!)^k} (b(1+\log b))^{k-1} \vol_n(X).
$$
This proves Theorem \ref{TheoremE}.
\end{proof}

Note that for every $r>1-\frac{1}{k}$ the factor $(b(1+\log
b))^{\frac{k-1}{k}}$ grows less fast than $b^r$ when
$b\rightarrow\infty$.

If $\dim X=n=3p$ we can improve Theorem \ref{TheoremE} as follows:

\begin{theorem}\label{8.1}
Let $X$ be a compact, orientable
manifold, of dimension $n=3p$. Suppose there exist classes
$\beta_1, \beta_2, \beta_3\in H^p(X,\RR)$  such that
$\beta_1\cup\beta_2\cup\beta_3\not= 0$.   Then for every
Riemannian metric $g$ on $X$ we have 
$$
\stsys_p(g)\le\left({n\choose p}
\frac{n!}{(p!)^3}\right)^{\frac{1}{3}} b_p(X)^{\frac{2}{3}}
\vol_n(g)^{\frac{1}{3}}.
$$
\end{theorem}

\begin{proof}
Set $\lambda_1=\lambda_1(H_p(X,\ZZ)_{\rr}, \strichleer_{L^2})$.  Let
$b=b_p(X)$, and set also $\lambda^*_b=\lambda_b(H^p(X,\ZZ)_{\rr},
\strichleer_{L^2})$.  We choose linearly independent classes
$\alpha_1, \ldots, \alpha_b$ in $H^p(X,\ZZ)_{\rr}$ such that
$|\alpha_i|_{L^2}\le\lambda^*_b$.  As in the proof of Theorem
\ref{TheoremE} we choose indices $i_1, i_2, i_3$ in $\{1,\ldots,b\}$
such that
$$
\alpha_{i_1}\cup\alpha_{i_2}\cup\alpha_{i_3}\not= 0.
$$
We may assume that $X$ is connected and oriented.
We consider $h=\PD_{\rr} (\alpha_{i_1}\cup\alpha_{i_2})\in
H_p(X,\ZZ)_{\rr} \setminus\{0\}$, harmonic forms $\omega_{i_1},
\omega_{i_2}$ representing $\alpha_{i_1},\alpha_{i_2}$ and a closed
$p$-form $\omega$ representing an arbitrary $\alpha\in H^p(X,\RR)$.
Then we have
$$
[h,\alpha]=\int_X \omega_{i_1} \wedge\omega_{i_2}
\wedge\omega\le\frac{n!}{(p!)^3} \parallelom_{\infty}\int_X
|\omega_{i_1}||\omega_{i_2}| d\vol_n,
$$
cf.~\cite{[Fe1]}, 1.7.5 and 1.8.1 for the constant.  Using the
definition of comass $\parallelh=\sup\limits_{\parallelalf^*\le
1}[h,\alpha]$, we obtain
$$
\parallelh\le\frac{n!}{(p!)^3} |\alpha_{i_1}|_{L^2}
|\alpha_{i_2}|_{L^2} \le \frac{n!}{(p!)^3} (\lambda^*_b)^2.
$$
We apply (\ref{(5.5)}) to the fundamental cycle $z$ of $M$ to obtain
$$
\stsys_p(g)\le{n\choose p}^{\frac{1}{2}} \lambda_1 \vol_n(g)^{\frac{1}{2}}.
$$
The preceding inequalities imply
$$
\stsys_p(g)^3 \le \stsys_p(g)^2 \parallelh \le {n\choose p} \frac{n!}{(p!)^3}
(\lambda_1 \lambda^*_b)^2 \vol_n(g).
$$
Now we apply Banaszczyk's inequality $\lambda_1 \lambda^*_b \le b$,
cf.~(\ref{(7.1)}), to complete the proof.
\end{proof}


\appendix                               






\ack
The authors have benefited from helpful
discussions with M.~Kreck and express appreciation to L.~Ambrosio,
G.~Dula, J.~Lagarias, C.~LeBrun, F.~Morgan, S.~Weinberger and B.~White
for insightful comments.



\frenchspacing
\bibliographystyle{plain}





\end{document}